\documentclass[12pt]{article}

\usepackage{amsmath}
\usepackage{amssymb}
\usepackage{amsfonts}
\usepackage{amsthm}
\usepackage{latexsym}
\usepackage{epsfig}
\usepackage{array}
\usepackage{boxedminipage}

\newcommand{\R}{\mathbb{R}}

\newtheorem{theor}{Theorem}
\newtheorem{lema}{Lemma}
\newtheorem{propo}{Proposition}

\newcommand{\vs}{\vspace{1mm}}
\newcommand{\vv}{\vspace{2mm}}
\newcommand{\vvv}{\vspace{3mm}}
\newcommand{\vvvv}{\vspace{4mm}}
\newcommand{\vvvvv}{\vspace{5mm}}

\newcommand{\rk}{\mathrm{rk} \hspace{0.5mm}}
\newcommand{\cork}{\mathrm{cork} \hspace{0.5mm}}
\newcommand{\ke}{\mathrm{ker} \hspace{0.5mm}}
\newcommand{\im}{\mathrm{im} \hspace{0.5mm}}
\newcommand{\re}{\mathrm{Re} \hspace{0.5mm}}

\newcommand{\tr}{{\sf T}}

\newcommand{\op}{\overline{p}}
\newcommand{\ou}{\overline{u}}

\newcommand{\oq}{\overline{q}}
\newcommand{\ow}{\overline{w}}

\newcommand{\dsp}{\displaystyle}

\begin{document}

\begin{center}
{\Large\bf Transcritical bifurcation without 
parameters \vspace{2mm} \\ 
in memristive circuits\footnote{This work was supported by Research Project 
MTM2015-67396-P of Ministerio de Econom\'{\i}a y Competitividad (MINECO)/Fondo 
Europeo de Desarrollo Regional (FEDER).}}

\vvv

{\sc 
Ricardo Riaza}\\ 
\ \vspace{-3mm} \\
Depto.\ de 
Matem\'{a}tica Aplicada 
a las TIC - 
ETSI 
Telecomunicaci\'{o}n  \\ 
Universidad Polit\'{e}cnica de Madrid,  
Spain \\
{\sl ricardo.riaza@upm.es} \\
\end{center}

\begin{abstract}
The transcritical bifurcation without parameters (TBWP) describes a 
stability change along a line of equilibria, resulting from the 
loss of normal hyperbolicity at a given point of such a line. 
Memristive circuits systematically yield 
manifolds of non-isolated equilibria, and in this paper
we address a systematic characterization of the TBWP 
in circuits with a single memristor. To achieve this
we develop two mathematical results of independent interest;
the first one is an extension of the TBWP theorem to explicit
ordinary differential equations (ODEs) in arbitrary dimension;
the second result drives the characterization of this
phenomenon to semiexplicit differential-algebraic
equations (DAEs), which provide the appropriate framework for
the analysis of circuit dynamics.
In the circuit context the analysis is performed in graph-theoretic
terms: in this setting, our first working scenario is restricted
to passive problems (exception made of the bifurcating memristor),
and in a second step some results are presented 
for the analysis of non-passive cases. The latter context
is illustrated by means of a
memristive neural network model.
\end{abstract}



\vvvv



\vvvv

\section{Introduction}
\label{sec-intr}

Memristors and other related electronic 
devices \cite{chuamemristor71,diventra09,memristor2008} 
are known to exhibit systematically manifolds of non-isolated equilibrium
points. This is a consequence of the form of the voltage-current characteristic of memristive
devices. As shown in \cite{stabmem}, in the absence of certain configurations
equilibrium manifolds of strictly locally passive memristive circuits are normally hyperbolic,
that is, all remaining eigenvalues of the linearized vector field (except
for those whose
eigenvectors span the tangent space to the equilibrium manifold) are away from the imaginary
axis. From a qualitative point of view it is therefore natural to examine what happens 
in the memristive circuit dynamics when the aforementioned passivity assumption does not hold.

This problem must be framed in the theory of {\em bifurcation without parameters} originally
introduced in the seminal papers \cite{fiedler00a, fiedler00b, fiedler00c}; cf.\ also the recent
book \cite{liebscherbook}. When normal hyperbolicity fails, a change in the local qualitative
properties typically occurs along the equilibrium manifold, hence the 
``bifurcation without parameters'' term. In this context, the most basic 
phenomenon is the transcritical bifurcation without parameters (TBWP throughout the whole manuscript),
which describes the transition
of one eigenvalue through the origin under certain local conditions on the vector field. Our purpose
in this paper is to present a systematic circuit-theoretic
characterization of this bifurcation for memristive
circuits. Due to the systematic presence of non-isolated
equilibria in memristive circuits, this 
is the most elementary phenomenon
responsible for a stability loss in nonlinear circuits with memristors,
and its analysis seems therefore to be very relevant for the development
of the qualitative theory of memristive circuits.

This is actually a three-fold goal. First, the characterization of the TBWP
in \cite{fiedler00a} (and also in \cite{liebscherbook}) is only addressed for two-dimensional dynamics.
However, most nonlinear circuits involve a large number of dynamic variables and a
two-dimensional model reduction is rarely feasible.
For their results to apply to higher
dimensional problems, the authors assume in \cite{fiedler00a, liebscherbook}  
that a prior reduction to a center two-dimensional
manifold has been performed, but no explicit
conditions paving the way
for an appropriate reduction are given in arbitrary (finite) dimension. 
As a somewhat natural (yet not trivial)
extension of their
characterization we will present a TBWP theorem for explicit ODEs in $\R^n$, addressing the
geometrical conditions which allow for a center manifold reduction in which the two-dimensional
conditions of \cite{fiedler00a} do hold. 
This first goal is addressed in Section \ref{sec-Rn} (cf.\ Theorem
\ref{th-Rn}).

Many nonlinear circuits do
not admit a dynamical description in terms of an
explicit ODE. This is clearly the case in large scale integration circuits, for which
such a state space description in terms of an explicit ODE model is hardly automatable.
For this reason, semistate models based on differential-algebraic equations (DAEs) are preferred
instead \cite{et00, gun1, LMTbook, wsbook, carentop}. Analytical results involving dynamical
systems intended to apply to nonlinear circuit models should therefore be addressed to DAEs.
Many qualitative investigations about nonlinear circuits
require a prior reduction
to an ODE model, involving unnecessarily restrictive hypotheses which, in addition,
make the analysis more difficult (see e.g.\ \cite{chua80} as a sample). 
In the DAE framework the approach is different:
instead of trying to 
drive a model to the ODE context in order to apply a given known result, 
it is more convenient to extend such result to the DAE setting, allowing for a direct
application to whatever semistate model. In this direction, 
our second goal is to drive the TBWP theorem to the semiexplicit DAE setting, a task which is 
accomplished in Section \ref{sec-daes} and,
specifically, in Theorem \ref{th-daes}. 

As indicated above, our third goal is to obtain a characterization 
of the TBWP for memristive circuits
in circuit-theoretic terms. This means that the characterization should be stated in terms
of the underlying circuit digraph and the electrical features of the devices.
This so-called {\em structural} approach has its roots in the state-space formulation problem
(whose origins can be traced back to \cite{bashkow57, bryant62}) and, more recently, has been
successfully applied to the DAE index characterization of several nonlinear 
circuit  models \cite{et00, gun1, iwataMC, wsbook, takamatsu2014, carentop, carenhabil}.
This approach makes it possible to directly transfer different analytical
and qualitative results to circuit simulation 
programs. Allowed by the
TBWP theorem for DAEs obtained in Section \ref{sec-daes}, 
such a characterization is detailed in Section \ref{sec-circuits} for circuits
displaying a line of equilibria (that is, including exactly one memristor), under
the assumption that the failing of a passivity assumption on this memristor is the
one responsible for the loss of normal hyperbolicity: cf.\ Theorem
\ref{th-circuits}. The analysis in this section extends some preliminary results
presented in \cite{aims2014}.

Section \ref{sec-nonpassive} discusses this phenomenon relaxing the
passivity assumption on the remaining circuit devices, and includes an example
coming from the theory of memristive neural networks.
Finally, Section \ref{sec-con}
compiles some concluding remarks.


\section{The TBWP theorem for explicit ODEs}
\label{sec-Rn}

\subsection{Two-dimensional dynamics}

We begin by recalling the characterization of the TBWP in two-dimensional problems
presented by Fiedler, Liebscher and Alexander in \cite{fiedler00a}.

\begin{theor}[Fiedler, Liebscher \& Alexander, 2000] \label{th-fiedler}
Consider the system
\begin{subequations} \label{fiedler}
\begin{eqnarray}
x' & = & \xi_1(x,y) \\
y' & = & \xi_2(x,y),
\end{eqnarray}
\end{subequations}
with $\xi \in C^2(\R^2, \R^2)$, and assume that
\begin{enumerate}
\item $\xi(x, 0)=0$;
\item $\dsp\frac{\partial \xi_2}{\partial y}(0,0) = 0$; 
\item $\dsp\frac{\partial \xi_1}{\partial y}(0,0) \neq 0$;
\item $\dsp\frac{\partial^2 \xi_2}{\partial x \partial y}(0,0) \neq 0$.
\end{enumerate}
Then (\ref{fiedler}) is locally orbitally $C^1$-equivalent 
to the normal form
\begin{subequations} \label{eqfn}
\begin{eqnarray}
x' & = & y \\
y' & = & xy
\end{eqnarray}
\end{subequations}
around the origin.
\end{theor}

\vspace{2mm}

\noindent Needless to say, it is enough to assume
that condition 1 holds for $x$ sufficiently close to 0.
Note that $y=0$ is a line of equilibria for both (\ref{fiedler}) and 
(\ref{eqfn}), and that $\lambda=0$ is an eigenvalue for the linearization
of both systems at any $(x,0)$. This zero eigenvalue becomes a double (index-two) one
at $(0,0)$, in a way such that the second system eigenvalue changes sign along
the line of equilibria; specifically, this second eigenvalue is positive (resp.\ negative)
if $x>0$ (resp.\ $x<0$) in the normal form (\ref{eqfn}). 
This means that the line of equilibria is normally hyperbolic 
for $x \neq 0$, and  a stability change along the line
of equilibria occurs as a result of the loss of normal hyperbolicity
 at the origin. This is the transcritical bifurcation without
parameters.

\subsection{The TBWP for explicit ODEs in $\R^n$}

Theorem \ref{th-fiedler} can be extended to explicit ODEs in 
arbitrary (finite) dimension as follows. Mind a notational abuse
used in different situations throughout the paper, namely writing both
$f(x)$ and $f(x_1, \ldots, x_n)$, the latter standing for the (more cumbersome)
$f((x_1, \ldots, x_n))$. Obviously, the result below can be also stated for an open set $\Omega \subseteq 
\R^n$ with $0 \in \Omega$ or for a germ of a map at the origin.

\vs

\begin{theor}[TBWP in $\R^n$] \label{th-Rn}
Assume that $f \in C^2(\R^n, \R^n)$ verifies the following.
\begin{enumerate}
\item {$f(x_1, 0, \ldots, 0)=0$}; 
\item $f'(0)$ has a double index-two zero eigenvalue, and {$\re \lambda \neq 0$}
for the remaining ones;
 
\item {$f''(0)pq \notin \im f'(0)$}, if $p \in \ke f'(0)-\{0\}$ and $q \in \ke (\hspace{0.5mm}f'(0))^2-\ke f'(0).$
\end{enumerate}
Then there exists a local, two-dimensional, $C^2$ center manifold,
where the reduced dynamics admits a description in local coordinates
of the form $u'=\xi(u)$, with $\xi$ verifying the conditions of Theorem \ref{th-fiedler} in
$u^*=(0,0)$.
\end{theor}

\vspace{2mm}

\noindent {\bf Proof.} The proof relies on the fact that the linear transformation
driving the linear part to Jordan form leaves the line of equilibria invariant;
note also that condition 3 captures the geometric (transversality) hypothesis
which extends to higher-dimensional contexts the second-order condition in Theorem \ref{th-fiedler}.
For better clarity we proceed in numbered steps.

\vs

\noindent 1. Let $P$ be a matrix that drives $f'(0)$ to Jordan form:
\begin{equation*}
\hspace{1.3mm}\tilde{\hspace{-1.3mm}J}=P^{-1}f'(0)P = \begin{pmatrix} J_0 & 0 \\ 0 & J_h\end{pmatrix}, \ J_0 = \begin{pmatrix} 0 & 1 \\ 0 & 0\end{pmatrix}, \ J_h \text{ hyperbolic.}
\end{equation*}

\vs

\noindent 2. Under the change of coordinates $x=Py$, the system
$x'=f(x)$ is transformed into $y'=\hspace{1.0mm}\tilde{\hspace{-1.0mm}f}(y)=P^{-1}f(Py)$,
which reads as
\begin{eqnarray*}
u' & = & \hspace{1.0mm}\tilde{\hspace{-1.0mm}f}_{\hspace{-0.3mm}1}(u, v) = J_0 u + \eta(u, v) \\
v' & = & \hspace{1.0mm}\tilde{\hspace{-1.0mm}f}_{\hspace{-0.3mm}2}(u, v) = J_h v + \gamma(u, v)
\end{eqnarray*}
with $y=(u, v), \ u \in \R^2, \ v \in \R^{n-2}$ and $\eta'(0,0)=0$, $\gamma'(0,0)=0$.

\vs

\noindent 3. Now, from the fact that 
$\hspace{0.6mm}\tilde{\hspace{-0.6mm}p} \in \ke
\hspace{1.0mm}\tilde{\hspace{-1.0mm}f}'(0) 
 \Leftrightarrow p=P\hspace{0.8mm}\tilde{\hspace{-0.8mm}p} 
\in \ke f'(0)$ it follows that $P$ leaves the space $\ke f'(0)=$ span\hspace{0.5mm}$\{e_1\}$
(with $e_1=(1, 0, \ldots, 0)$)
invariant, and therefore the equilibrium line of
$\hspace{1.0mm}\tilde{\hspace{-1.0mm}f}$ is $(u_1, 0, \ldots, 0)$.

\vs

\noindent  4. Using, analogously, the properties
$\hspace{0.7mm}\tilde{\hspace{-0.7mm}q} \in \ke
(\hspace{0.7mm}\tilde{\hspace{-0.7mm}f}'(0))^2
 \Leftrightarrow q=P\hspace{1.0mm}\tilde{\hspace{-1.0mm}q} 
\in \ke (\hspace{0.2mm}f'(0))^2$
and $\hspace{1.0mm}\tilde{\hspace{-1.0mm}w} \in \im
\hspace{1.0mm}\tilde{\hspace{-1.0mm}f} '(0) 
 \Leftrightarrow w=P\hspace{1.0mm}\tilde{\hspace{-1.0mm}w} 
\in \im f'(0)$,
the condition $$f''(0)pq \notin \im f'(0) \text{ with }p \in \ke f'(0)-\{0\} \text{ and }
q \in \ke (\hspace{0.5mm}f'(0))^2-\ke f'(0)$$ 
yields $$\hspace{1.0mm}\tilde{\hspace{-1.0mm}f}''(0)\hspace{1.0mm}\tilde{\hspace{-0.7mm}p}\hspace{0.7mm}\tilde{\hspace{-0.7mm}q} \notin \im \hspace{1.0mm}\tilde{\hspace{-1.0mm}f}'(0), \text{ with } \hspace{0.7mm}\tilde{\hspace{-0.7mm}p} \in \ke \hspace{1.0mm}\tilde{\hspace{-1.0mm}f}'(0)-\{0\}
\text{ and } \hspace{0.7mm}\tilde{\hspace{-0.7mm}q} \in \ke (\hspace{0.5mm}\hspace{1.0mm}\tilde{\hspace{-1.0mm}f}'(0))^2-\ke \hspace{1.0mm}\tilde{\hspace{-1.0mm}f}'(0)$$
and, in turn, this leads to
\begin{equation} \label{segundo}
\dsp\frac{\partial^2 \eta_2}{\partial u_1 \partial u_2}(0,0) \neq 0.
\end{equation}

\vv

\noindent 5. The system $y'=\hspace{1.0mm}\tilde{\hspace{-1.0mm}f}(y)$ admits a local center manifold
of the form $v=\zeta(u)$, 
with $\zeta(0)=0, \ \zeta'(0)=0$ (see e.g.\ \cite{carr, perko, wiggins}). 
The dynamics on this manifold reads as
$$ u' = \xi (u) = J_0 u + \eta(u, \zeta(u)).$$

\vv

\noindent 6. Locally, the curve of equilibria must belong to the center manifold \cite{carr,wiggins},
and this yields condition 1 of Theorem \ref{th-fiedler}, that is,
$\xi(u_1, 0)=0$. Additionally, the form of
$J_0$ renders conditions $2$-$3$ trivial,
and (\ref{segundo}) yields condition 4, that is,
$$\dsp\frac{\partial^2 \xi_2}{\partial u_1 \partial u_2}(0,0) \neq 0,$$
since $\eta'(0,0)=0$, $\zeta'(0,0)=0$. This completes the proof.

\hfill $\Box$

\vs

\noindent {\bf Remark.} Geometrically, condition 3 expresses the transversality at $x^*$
of the center manifold and the
so-called singular manifold 
\begin{equation} \label{singular}
\{x \in \R^n \ / \ \det f'(x) = 0\},
\end{equation}
 as a consequence
of the following well-known property from matrix analysis (we omit the proof;
find details e.g.\ in \cite{qsib}).

\begin{lema}\label{lema-det} 
If $H \in C^1(\R^m, \R^{n \times n})$, $\rk H(x^*) = n-1$ and $p \in \ke H(x^*)-\{0\}$, then
$$ (H'(x^*)q)p \notin \im H(x^*) \Leftrightarrow (\det H)'(x^*)q \neq 0.$$
\end{lema}

\vs

\noindent Lemma \ref{lema-det} implies in particular that, if
$f\in C^2(\R^n, \R^n)$ and $\rk f'(x^*) = n-1$, then
$$f''(x^*)pq \notin
\im f'(x^*) \Leftrightarrow (\det f')'(x^*)q \neq 0,$$
a condition which, in the setting of Theorem \ref{th-Rn}, 
expresses the transversal intersection of the direction spanned by the
generalized eigenvector $q$ (hence
of the center manifold itself) and the singular manifold (\ref{singular}),
as indicated in the Remark above. For later use,
we also note that for a single-parameter valued matrix map
$H \in C^1(\R, \R^{n \times n})$ with $H'(\lambda^*)=n-1$ we have
$$H'(\lambda^*)p \notin \im H(\lambda^*) \Leftrightarrow (\det H)'(\lambda^*) \neq 0.$$

Lemma \ref{lema-det} also shows that condition 3 does not depend
on the choice of $q$, since
$\ke f'(0)$ is tangent to
the singular manifold, so that 
$f''(0)pp \in \im f'(0)$
and then, for $\hat{q}=\alpha q + \beta p$
with $\alpha \neq 0$, 
we have $${f''(0)p\hat{q}} = (\alpha f''(0)p q + \beta f''(0)p p)
\notin \im f'(0) \Leftrightarrow f''(0)p q \notin \im f'(0).$$

Note also that the form 
of $\ke f'(0)=$ span$\{e_1\}$ makes it possible to simplify the statement of condition
3 to
$${f_{x_1 x}(0)q \notin \im f'(0)},$$
where $f_{x_1 x}$ denotes the matrix of partial
derivatives $\dsp\left( \frac{\partial^2 f_i}{\partial x_1 \partial x_j}\right)$.

Finally, from the \u{S}o\u{s}ita\u{\i}\u{s}vili-Palmer Theorem \cite{palmer},
it follows that the normal form for the TBWP in $\R^n$ is
\begin{eqnarray*}
x' & = & y \\
y' & = & xy  \\
v' & = & J_h v. 
\end{eqnarray*}
Certainly, the form of the latter equation may be further simplified to that of
a standard node or saddle point, depending on the inertia of $J_h$.

\section{TBWP in semiexplicit DAEs}
\label{sec-daes}

Along the route indicated in Section \ref{sec-intr}, we extend below the TBWP to
the setting of semiexplicit DAEs.

\begin{theor}[TBWP in semiexplicit index-one DAEs]\label{th-daes}
Let $h \in C^2(\R^{r+p},\R^r)$, $g \in C^2(\R^{r+p},\R^p)$, and consider the system
\begin{subequations} \label{dae}
\begin{eqnarray}
y' & =&  h(y, z) \\
0 & = & g(y, z). \label{daeb}
\end{eqnarray}
\end{subequations}
Write $E=
\left(\begin{array}{cc}I_r & 0 \\ 0 & 0\end{array}\right)$, \ $F=(h,g)$.
Assume that $g_z(0,0)$ is non-singular and that
\begin{enumerate}
\item {$h(y_1, 0, 0)=0$, $g(y_1, 0, 0)=0$}; 

\item the matrix pencil $\lambda E - F'(0,0)$ has a double index-two zero eigenvalue, 
and {$\mathrm{Re} \hspace{0.5mm} \lambda \neq 0$} for the remaining eigenvalues;

\item {$F''(0,0)\op\hspace{0.1mm}\oq \notin \im F'(0,0)$}, where
\begin{equation}
\op \in \ke F'(0, 0)-\{0\}, \ 
\oq \in \ke (F'(0,0))^2-\ke F'(0,0).\label{opoq} \end{equation}
\end{enumerate}
Then, there exists an invariant, two-dimensional, $C^2$ submanifold of
$g(y,z)=0$ where the dynamics admits a local description of the form
$u'=\xi(u)$ with
$\xi$ satisfying the conditions of Theorem \ref{th-fiedler} at the origin.
\end{theor}

\vspace{2mm}

\noindent Before proceeding with the proof we present some auxiliary results.

\begin{lema}[Schur] \label{lema-schur}
Let $D$ be a non-singular matrix and
\begin{equation}M=\begin{pmatrix}A & B \\ C & D\end{pmatrix}, \ (M/D)=A-BD^{-1}C,\label{schur}
\end{equation}
with $A$ (hence $M$) square. Then $\det M = \det (M/D) \det D$ and $\cork M = \cork (M/D).$
\end{lema}

We will make use of this Lemma at several points in our analysis, mostly with 
\begin{eqnarray} \label{MFprima}
 M = F'(0,0)= \left(\begin{array}{cc}h_y(0,0) & h_z(0,0) \\ g_y(0,0) & g_z(0,0) \end{array}\right)
\equiv \left(\begin{array}{cc}A & B \\ C & D \end{array}\right).\end{eqnarray}

The proof of Theorem \ref{th-daes} will be based on checking conditions 1-3 of Theorem \ref{th-Rn}
for the reduced dynamics of the DAE (\ref{dae}) on the solution manifold (\ref{daeb}). Conditions
1 and 2 will be derived in a more or less straightforward manner; condition 3 is not trivial, though.
Remember that the goal is to state the conditions in terms of the original
problem setting, that is, in terms of $h$ and $g$ (that is, of $F$), as it is done in our
statement of condition 3. 
But note  that it is the matrix pencil $\lambda E - F'(0,0)$ (and not the matrix $F'(0,0)$) 
the one that is assumed to
have a double, index-two zero eigenvalue; this means that it is not obvious that there
should exist an 
$\oq$ satisfying the requirement depicted in (\ref{opoq}) in light of the previous hypotheses.
As a cautionary example consider, for instance, the Schur reduction of
$$M = \begin{pmatrix} 0 & 1 & 1 \\ 1 & 0 & 0 \\ 0 & 1 & 1 \end{pmatrix}, \ \
 (M/D) = \begin{pmatrix} 0 & 0 \\ 1 & 0 \end{pmatrix},$$ 
and note that $\lambda=0$ is a double eigenvalue for $(M/D)$ but a simple one for $M$; no
generalized eigenvector exists in this case for $M$. 
This cannot occur in the setting of Theorem \ref{th-daes}
(that is, there will indeed exist a $\oq$ satisfying the condition in (\ref{opoq}))
because of item (c) of Lemma \ref{lema-aux} below.

\begin{lema} \label{lema-aux}
Given $M$ and $(M/D)$ in (\ref{schur}), consider the 
operators $L$ and 
$T: \R^r \to \R^{r+p}$ defined by
$$Lu = \begin{pmatrix}u \\ -D^{-1}Cu \end{pmatrix},
\ \ Tu=\begin{pmatrix}u \\ 0 \end{pmatrix}.$$
Then

\noindent (a) $\op \in \ke M \Leftrightarrow \op = Lp,$ with $p \in \ke (M/D);$

\vs

\noindent (b) $w \in \im (M/D) \Leftrightarrow Tw 
\in\hspace{-0.2mm} \im M$. Actually,
$Tw=M\ou \Leftrightarrow \ou = Lu$, with $w = (M/D)u$;

\vs

\noindent (c) if $\ke M  \subseteq \im T= \R^r\times \{0\}$, 
then
$\overline{q}\in \ke M^2 \Leftrightarrow \oq = Lq$, with $q \in \ke (M/D)^2$.
\end{lema}

\vv

\noindent {\bf Proof.} Items (a) and (b) are immediate in light of the definition of
$(M/D)$ in (\ref{schur}). Regarding (c), let {$\overline{q}\in \ke M^2$}
and denote $\ow = M\hspace{0.1mm}\oq \in \ke M \cap \im M$. 
Then:

\noindent  (i) since $\ow \in \ke M$, then $\ow = Lw$
with $w \in \ke (M/D)$ because of (a);

\noindent (ii) owing to the hypothesis $\ke M  \subseteq \im T$,
necessarily $\ow=Tu$ for a certain $u$; additionally,
because of the form of $L$ and $T$ it follows that $u=w$ and then
$\ow=Lw=Tw$;

\noindent  (iii) the condition $\ow = T w = M\hspace{0.2mm}\oq$ implies, because of (b), 
that {$\oq = Lq$}, with $w = (M/D)q$;

\noindent  (iv) finally, since $w \in \ke (M/D)$ (cf.\ (i)), it follows that 
{$q \in \ke (M/D)^2$}. 

The converse result in (c) is entirely analogous and the proof is left to the reader.

\hfill $\Box$

\vv

\noindent
With $M$ as defined in (\ref{MFprima}), then the linear operator
$L$ is the differential at the origin of the parameterization 
$y \to (y, \psi(y))$ of the manifold ${\cal M}$ in (\ref{daeb}), with $\psi$ given
by the implicit function theorem, and therefore $L$
will define an isomorphism
$\R^r \to \ke \begin{pmatrix} C & \hspace{-1.5mm} D \end{pmatrix} = T_{(0,0)} {\cal M}$. 
Item (a) in Lemma \ref{lema-aux} expresses that $L$ also induces an isomorphism 
between the spaces $\ke (M/D) \to \ke M \subseteq T_{(0,0)} {\cal M}$.
Moreover, in the scenario assumed in
(c), one can check that $\ke M^2 \subseteq T_{(0,0)} {\cal M}$ and
that $L$ also induces an isomorphism $\ke (M/D)^2 \to \ke M^2$.

For later use, we also note that a coordinate description $\alpha$ of a map $\beta$ 
defined on ${\cal M}$, that is, a relation of the form
$\alpha(y)=\beta(y, \psi(y))$, implies $\alpha'(0) =\beta'(0,\psi(0)) L$. We will make use of
this remark in the final step of the proof of Theorem \ref{th-daes}.

\vv

\noindent {\bf Proof of Theorem \ref{th-daes}.}

\vs

\noindent 1. The hypothesis that $g_z(0,0)$ is non-singular implies, by the implicit function 
theorem, that 
${\mathcal M} \equiv g=0$ is locally a manifold that can be described by
$z=\psi(y)$.  The goal is then to apply Theorem \ref{th-Rn} to the reduced system
\begin{equation}y'=f(y)=h(y, \psi(y)).\label{reduced}\end{equation}
Specifically, we need to check that the requirements imposed on $F=(h, g)$ 
yield the conditions 1-3 in Theorem \ref{th-Rn}.

 \vs

\noindent 2. The first condition holds trivially, since
$$f(y_1, 0)=h(y_1, 0, \psi(y_1, 0)) = h(y_1, 0, 0) = 0,$$ where the second identity
is due to the fact that $\psi(y_1, 0)=0$ because $g(y_1, 0, 0)=0$.

\vs

\noindent 3. The second condition also follows easily from the implicit function
theorem, since the linearization of the reduced system (\ref{reduced}) at the origin is
$$f'(0)=h_y(0,0)+h_z(0,0)\psi'(0) = h_y(0,0)- h_z(0,0)(g_z(0, 0))^{-1}g_y(0,0),$$
and the spectrum of $f'(0)=A-BD^{-1}C$ equals that of the matrix pencil $\lambda E-M$, as an immediate
consequence of Schur's lemma:
\begin{eqnarray*}
&& \det (\lambda E-M)   =  \det \left( \lambda
\left(\begin{array}{cc}I_r & 0 \\ 0 & 0\end{array}\right) -
\left(\begin{array}{cc}A & B \\ C & D \end{array}\right)\right)\\ 
 && \hspace{3mm} = \det 
\left(\begin{array}{cc}\lambda I_r -A  & -B \\ -C & -D \end{array}\right)
  = \det \left(\lambda I_r - (A- BD^{-1}C)\right) \det (-D).
\end{eqnarray*}

\vs

\noindent 4. The only non-immediate step consists in checking that
$$F''(0,0)\op\hspace{0.1mm}\oq \notin \im F'(0,0),$$ 
where
$\op \in \ke F'(0, 0)-\{0\}, \ 
{\oq \in \ke (F'(0,0))^2-\ke F'(0,0)},$ implies
$$f''(0)pq \notin \im f'(0),$$ 
with $p \in \ke f'(0)-\{0\}$ and ${q \in \ke (\hspace{0.5mm}f'(0))^2-\ke f'(0)}.$
Note that item (c) of Lemma \ref{lema-aux} applies because
$\ke F'(0,0)=$ span$\{e_1\}$, and then
$\oq$ yields a generalized eigenvector $q$ of $f'(0)$, with $\oq=Lq$.

It then suffices to use the characterization
$${F''(0,0)\op\hspace{0.1mm}\oq \notin \im F'(0,0) \Leftrightarrow (\det F')'(0,0)\oq \neq 0},$$
following from Lemma \ref{lema-det} because, differentiating
$\det F'=\det g_z \det (F'/g_z)$ and using
the fact that 
$\det (F'/g_z)(0,0)=\det (\hspace{0.5mm}f')(0)=0$ (because $\lambda=0$ is an eigenvalue), 
we have
$$(\det F')'(0,0)=\det g_z(0,0)(\det (F'/g_z))'(0,0);$$ 
additionally, since $\det g_z(0,0) \neq 0$, it follows that
$$ { (\det F')'(0,0)\oq \neq 0} \Leftrightarrow (\det (F'/g_z))'(0,0)Lq \neq 0 \  
{ \Leftrightarrow  (\det (\hspace{0.5mm}f'))'(0)q \neq 0},$$
because
$\det \hspace{0.5mm}f'(y)=\det (F'/g_z)(y, \psi(y))$
$\Rightarrow (\det (\hspace{0.5mm}f'))'(0)=(\det (F'/g_z))'(0,0)L$ as indicated
right before the proof of Theorem \ref{th-daes}. Finally,
\begin{eqnarray*}
{(\det (\hspace{0.5mm}f'))'(0)q \neq 0 \Leftrightarrow f''(0)pq \notin \im f'(0)} 
\end{eqnarray*}
and the proof is complete.

\mbox{} \ \hfill $\Box$

\vv

\noindent {\bf Remark.} As in the explicit ODE case,
condition 3 can be recast as
$${F_{x_1 x}(0, 0)\oq \notin \im F'(0, 0)},$$
where $x=(y,z)$ and
$F_{x_1 x}$ stands for the matrix of partial derivatives 
$\dsp\left( \frac{\partial^2 F_i}{\partial x_1 \partial x_j}\right)$. 

\section{TBWP in memristive circuits}
\label{sec-circuits}

The memristor (an abbreviation for {\em memory-resistor})
is a new electronic device governed by a nonlinear, $C^1$ flux-charge relation  of the form
$\varphi = \phi(q)$. The existence of such a device was 
predicted for symmetry reasons by Chua in 1971 \cite{chuamemristor71}, but it was
not until
 2008 that it began to attract considerable attention. The reason for this
was the report in \cite{memristor2008} of the design of a nanometric memristor by the HP company.
The key aspect of this device is that, by differentiation of the aforementioned constitutive
relation, one gets the voltage-current relation 
\begin{equation} \label{memr}
v = M(q) i,
\end{equation}
with $M(q)=\phi'(q)$. For later use, we will assume that $\phi$ is a $C^2$ map.
Note that in (\ref{memr}) the ``resistance'' (or, better, {\em memristance})
$M$ depends on 
$q(t) = \int_{-\infty}^t i(\tau)d\tau$, so that the device somehow keeps track of its own
history (hence the memory-resistor name). A great amount of research has been directed
to this and other related devices since 2008; cf.\ \cite{adamatzky, diventra09, 
chuamemristor08, kavehei10, kim2012bis, 
messias10, muthus10, muthuschua10, pershin10c, memactive,
stabmem, tetzlaffBook} as a sample of literature.

It is easy to check  that memristors systematically yield manifolds
of non-isolated equilibrium points;
find details below. In order to focus on problems with {\em lines} of equilibria,
we will restrict our attention to circuits including a single memristor,
besides capacitors,
inductors, resistors, and 
(independent) voltage and current sources. Capacitors, inductors and resistors may be nonlinear, and
they will respectively be assumed to be defined by a voltage-dependent capacitance matrix $C(v_c)$,
a current-dependent inductance matrix $L(i_l)$, and a current-controlled $C^1$ description 
$v_r=\gamma(i_r)$ in the
case of resistors; for later use we denote the resistance matrix $\gamma'(i_r)$ as $R(i_r)$.
All three matrices need not be diagonal, allowing for the presence of coupling effects in the
corresponding sets of devices. Capacitors, inductors and resistors are said to be strictly locally
passive at a given operating point
if the corresponding characteristic matrix,
that is, $C(v_c)$, $L(i_l)$ or $R(i_r)$, 
is positive definite (a square matrix $P$
is positive definite if it verifies $v^{\tr} P v > 0$ 
for non-vanishing real vectors $v$; note that we do not require these matrices to be symmetric).  We assume
the circuit to be autonomous, that is, sources are DC ones or, in mathematical terms, they
take constant values (grouped together in vectors $V$ and $I$).

In this context, the circuit equations can be modeled by the differential-algebraic system
(see e.g.\ \cite{wsbook})
\begin{subequations}\label{branch}
\begin{eqnarray}
{q_m'} & {=} & {i_m} \label{branch0} \\
C(v_c) v_c'  & = & i_c \label{brancha} \\
L(i_l) i_l' & = & v_l \label{branchb} \\
0 & = & {B_m M(q_m)i_m} + B_c v_c + B_l v_l + B_r \gamma(i_r) + B_u V + B_j v_j  \label{branchc} \\
0 & = & {Q_m i_m} + Q_c i_c + Q_l i_l + Q_r i_r + Q_ui_u + Q_j I, \label{branchd}
\end{eqnarray}
\end{subequations}
where the subscripts $m$, $c$, $l$, $r$, $u$, $j$ are used for memristors, capacitors,
inductors, resistors, voltage sources and current sources, respectively. It is worth
emphasizing that (\ref{branchc}) and (\ref{branchd}) express Kirchhoff voltage and current
laws in terms of the so-called loop and cutset matrices $B$ and $Q$ (cf.\ the
Appendix and \cite{bollobas, 
wsbook, carentop}); we split these matrices as
$B=(B_m \ B_c \ B_l  \ B_r \ B_u \ B_j)$ and $Q=(Q_m \ Q_c \ Q_l  \ Q_r \ Q_u \ Q_j)$, where
$B_m$ (resp.\ $B_c,$ $B_l,$ $B_r$, $B_u$, $B_j$) corresponds to the columns of $B$
accommodating memristors (resp.\ capacitors, inductors, resistors,
voltage sources, current sources), and the same notational convention
applies to the cutset matrix.

By denoting $y=(q_m, v_c, i_l)$, \ $z=(i_m, i_c, v_l, i_r, v_j, i_u)$,
the DAE (\ref{branch}) takes the form
\begin{eqnarray*}
E(y)y' & =&  h(y, z) \\
0 & = & g(y, z)
\end{eqnarray*}
and equilibria are defined by the pair of conditions
$h(y,z)=0$, $g(y,z)=0$, that is, 
\begin{subequations} \label{eqconds}
\begin{eqnarray} 
i_m = i_c  =  v_l & = & 0 \\ 
B_c v_c + B_r \gamma(i_r)+ B_u V + B_j v_j & = & 0 \\
Q_l i_l + Q_r i_r+ Q_ui_u + Q_j I & = & 0. 
\end{eqnarray}
\end{subequations}
Note that the variable $q_m$ is not at all involved in
(\ref{eqconds}).
This means that, necessarily, no equilibrium point may be isolated,
since the variable $q_m$ unfolds any given equilibrium point
to a line (or even a higher dimensional set) of
equilibria. In our working setting equilibria will actually define
a line, as a consequence of the condition $\cork F'=1$ shown within the proof
of Theorem \ref{th-circuits}.

The previous remarks drive the stability analysis of equilibria in memristive
circuits to the mathematical context considered in \cite{aulbach, fiedler00a, 
fiedler00b, fiedler00c, liebscherbook}. 
In this setting, the existence of an $m$-dimensional manifold of equilibria 
implies that at least $m$ eigenvalues of the linearization
of the vector field at any of these equilibria are null.
The manifold is then said to be
{\em normally hyperbolic} (locally around such an equilibrium) 
if the remaining eigenvalues are not in the imaginary axis. 
The failing of the normal hyperbolicity requirement 
typically yields a bifurcation without parameters
\cite{fiedler00a, fiedler00b, fiedler00c, liebscherbook},
where the qualitative properties of the local phase portrait change.

In \cite{stabmem} one can find graph-theoretic 
conditions under which 
any manifold of equilibria of a strictly locally passive 
memristive circuit is guaranteed to be
normally hyperbolic. It is therefore natural
to address what happens when the passivity assumption does not hold;
allowed by Theorem \ref{th-daes}, 
Theorem \ref{th-circuits} below 
answers this question, again in circuit-theoretic
terms, for circuits with one memristor which becomes locally active at a given
operating point, yielding a transcritical bifurcation without parameters.
Note that in the statement of Theorem
\ref{th-circuits}, a VMC-loop is a loop composed only of voltage
sources, memristors and/or capacitors. ILC-cutsets, VML-loops, etc.\
are defined analogously.

\begin{theor}[TBWP in memristive circuits] \label{th-circuits}
Consider a nonlinear circuit with a single memristor, modeled by
(\ref{branch}). Fix an equilibrium point
$(q_m^*, v_c^*, i_l^*, i_m^*, i_c^*, v_l^*, i_r^*, v_j^*, i_u^*)$ (with 
$i_m^*=0,\ i_c^*=0, \ v_l^*=0$), and assume that the following conditions hold.
\begin{enumerate}
\item The circuit displays neither VMC-loops nor ILC-cutsets.

\item There is a unique VML-loop, which includes the memristor and at
least one inductor.

\item The capacitance, inductance and resistance matrices $C(v_c^*)$,
$L(i_l^*)$, $R(i_r^*)=\gamma'(i_r^*)$  are positive definite, with
$C(v_c^*)$ and $L(i_l^*)$ symmetric; additionally,
$M(q_m^*)=0$ and $M'(q_m^*) \neq 0$.
\end{enumerate}
Then the circuit undergoes a transcritical bifurcation without
parameters at the aforementioned equilibrium point; moreover,
near this bifurcating equilibrium,
all eigenvalues  of the linearization (but the null one)
have negative real part in the region where $M(q_m)>0$,
whereas a single (real) eigenvalue becomes positive at points
where $M(q_m)<0$.
\end{theor}

In the proof of Theorem \ref{th-circuits} we will make use of some
graph-theoretic results which are compiled in advance. 
Proofs of these auxiliary results can be found in \cite{andras1,
bollobas, wsbook, stabmem}.

\begin{lema}
\label{lema-siap} Let $B_i$ and $Q_i$ denote, for $i=1, \ 2, \ 3$, 
the submatrices
of $B$ and $Q$ defined by three pairwise-disjoint 
branch sets $K_1$, $K_2$, $K_3$ of a given directed graph. If $P$ is a 
positive definite matrix, then
\[\ke \left( \begin{array}{ccc} 
B_1 & 0 & B_3 P\\
0 & Q_2 & Q_3 
\end{array}
\right) = \ke B_1 \times \ke Q_2 \times \{0\}.\]
\end{lema}

The same terminological convention is used in Lemma \ref{lema-loopscutsets}
below. By a $K_i$-cutset (resp.\ loop) we mean a cutset (resp.\ loop)
defined only by branches belonging to $K_i$.
\begin{lema} \label{lema-loopscutsets}
The identity $\ke B_1 = \{0\}$ (resp.\ $\ke Q_2=\{0\}$) holds
if and only if the digraph has no {$K_1$-cutsets} 
(resp.\ {$K_2$-loops}).
\end{lema}

The proof of Theorem \ref{th-circuits} below is based
on Theorem \ref{th-daes}, which for simplicity 
was stated under the assumption
that the bifurcating equilibrium is located at the origin.
Obviously, we can make use of this result at a generic equilibrium $(y^*, z^*)$
(with $y^*=(q_m^*, v_c^*, i_l^*)$ and $z^*=(i_m^*, i_c^*, v_l^*, i_r^*, v_j^*, 
i_u^*)$) and we will do so without further explicit mention.

\vv

\noindent {\bf Proof of Theorem \ref{th-circuits}.} Note first that
the strict passivity assumption on $C(v_c^*)$ and $L(i_l^*)$ makes
these matrices non-singular, and therefore the maps $h$ and $g$ from
(\ref{dae}) have (at least locally) the form
\begin{subequations} \label{otra}
\begin{eqnarray*}
h(y, z) & = & \begin{pmatrix} i_m \\ (C(v_c))^{-1} i_c \\ (L(i_l))^{-1} v_l  \end{pmatrix}  \\
 \\
g(y,z) & = & \begin{pmatrix}
B_m M(q_m)i_m + B_c v_c + B_l v_l + B_r \gamma(i_r) + B_u V + B_j v_j \\
Q_m i_m + Q_c i_c + Q_l i_l + Q_r i_r + Q_ui_u + Q_j I
\end{pmatrix},
\end{eqnarray*}
\end{subequations}
where we are denoting $y=(q_m, v_c, i_l), \ z=(i_m, i_c, v_l, i_r, v_j, i_u)$.

\vs

\noindent 1. The matrix of partial derivatives $g_z(y^*, z^*)$ is 
(using $i_m^*=0, \ M(q_m^*)=0$)
$$g_z(y^*, z^*)=\begin{pmatrix}0 & 0 & B_l & B_r R(i_r^*)  & B_j& 0 \\ Q_m & Q_c & 0 & Q_r  & 0 & Q_u\end{pmatrix},$$
which is invertible in light of Lemmas \ref{lema-siap} and \ref{lema-loopscutsets}, since
$R(i_r^*)$  is positive definite and there are neither IL-cutsets 
(which are a particular instance of an ILC-cutset)
nor VMC-loops.

\vs

\noindent 2. Denoting $y=(q_m, \tilde{y})$, the conditions 
$h(q_m, \tilde{y}^*, z^*)=0$, $g(q_m, \tilde{y}^*, z^*)=0$
(arising in condition 1 of Theorem \ref{th-daes}) 
are trivially met (cf.\ (\ref{eqconds})).

\vs

\noindent 3. Condition 2 of Theorem \ref{th-daes} 
involves a matrix pencil spectrum which is given
by the determinant of
\begin{eqnarray}
\left( \begin{array}{ccccccccc}
 \lambda & 0 & 0 & -1 & 0 & 0 & 0 & 0 & 0 \vspace{0.3mm} \\ 
0 & \lambda I_c & 0 & 0 & \stackrel{}{-(C(v_c^*))^{-1}} & 0 & 0 & 0 & 0 \vspace{0.2mm} \\
0 & 0 & \lambda I_l & 0 & 0 & -(L(i_l^*))^{-1} & 0 & 0 & 0 \vspace{0.2mm}\\
0 &  -B_c &  0   &  0  & 0 & -B_l & -B_r R(i_r^*)  & -B_j& 0\vspace{0.2mm} \\ 
0 &   0  &  -Q_l & -Q_m & -Q_c & 0 & -Q_r  & 0 & -Q_u 
\end{array}\right), \label{bigmatrix}
\end{eqnarray}
and which 
can be written as
$\lambda d(\lambda)$ with $d(\lambda)=\det K(\lambda),$ provided
that 
 $K(\lambda)$ is the bottom-right submatrix of (\ref{bigmatrix})
(that is, the one
obtained after removing the top row and the left column).

The eigenvalue $\lambda=0$ being a double one amounts to 
$d(0)=0$ (ie.\ to $K(0)$ being singular) with $d'(0)\neq 0$.
Additionally, provided that ${\cork K(0)=1}$ (a condition which will be proved
to hold), then the condition $d'(0)\neq 0$ is equivalent (cf.\
Lemma \ref{lema-det}) to $K'(0)p \notin \im K(0)$
for $p \in \ke K(0)-\{0\}$. Finally, for the (double) zero eigenvalue
to be index-two it will be enough to show that
${\cork F'(y^*, z^*)=1}$. 
We examine this set of conditions in items 4, 5 and 6 below.

\vs

\noindent 4. The matrix $K(0)$ reads as
\begin{eqnarray*}
\hspace{-2mm}\begin{pmatrix}
0 & 0 & 0 & -(C(v_c^*))^{-1} & 0 & 0 & 0 & 0 \\
 0 & 0 & 0 & 0 & -(L(i_l^*))^{-1} & 0 & 0 & 0 \\
 -B_c &  0   &  0  & 0 & -B_l & -B_r R(i_r^*)  & -B_j& 0 \\ 
   0  &  -Q_l & -Q_m & -Q_c & 0 & -Q_r  & 0 & -Q_u
\end{pmatrix}
\end{eqnarray*}
and, following Schur's lemma (Lemma \ref{lema-schur}), this matrix is
easily seen to have the same corank as
\begin{eqnarray*}
\begin{pmatrix}
 B_c &  0   &  0  &  B_r R(i_r^*)  & B_j& 0 \\ 
   0  &  Q_l & Q_m &  Q_r  & 0 & Q_u
\end{pmatrix}.
\end{eqnarray*}

From  Lemmas \ref{lema-siap} and \ref{lema-loopscutsets}, the positive
definiteness of $R(i_r^*)$, the absence of IC-cutsets and the
existence of a unique VML-loop, it follows that
$\cork K(0)=1$, with
$$ \ke K(0) = \text{span} \{(0, p_l, p_m, 0, 0, 0, 0, p_u)\}$$
where $(p_l, p_m, p_u) \in \ke \begin{pmatrix} Q_l & \hspace{-1mm} Q_m & \hspace{-1mm} Q_u\end{pmatrix}$ (and note, for later use, that
${p_l \neq 0}$). 

\vs

\noindent 5. The condition $K'(0)p \notin \im K(0)$ is \vs
\begin{eqnarray*}
\begin{pmatrix}
0 \\ p_l \\ 0 \\ 0
\end{pmatrix} \notin \im\hspace{-1mm}
\begin{pmatrix}
 0 &  0 &  0 & -(C(v_c^*))^{-1} &  0 & \hspace{-1mm} 0 &  0 &  0 \\
  0 &  0 &  0 &  0 &  -(L(i_l^*))^{-1} &  0 &  0 &  0 \\
  -B_c &   0   &   0  &  0 &  -B_l &  -B_r R(i_r^*)  &  -B_j&  0 \\ 
    0  &   -Q_l &  -Q_m &  -Q_c &  0 &  -Q_r  &  0 &  Q_u
\end{pmatrix}.
\end{eqnarray*}
\vs Assuming that this condition does not hold, we would have a solution
for
\begin{eqnarray*}
B_c u_1 - B_l L(i_l^*) p_l + B_r R(i_r^*) u_6 + B_j u_7=0. 
\end{eqnarray*}
\vs
Together with the identity
 $Q_l p_l + Q_m p_m+ Q_up_u=0$ and the orthogonality of the
so-called cut and cycle spaces
$\ke B$, $\ke Q$ (cf. \cite{bollobas}), we would get
$p_l^{\tr}L(i_l^*)p_l=0$ and therefore 
$p_l=0$ (since $L(i_l^*)$ is positive definite), against
the fact that $p_l \neq 0$.

Hence, {$K'(0)p \notin \im K(0)$ for $p \in \ke K(0)-\{0\}$}, 
and as indicated above this implies that $\lambda=0$
is a double eigenvalue.

 \vs

\noindent 6. The matrix $F'(y^*, z^*)$ is
\begin{eqnarray*}
F'(y^*, z^*)=\begin{pmatrix}
0 & 0 & 0 & 1 & 0 & 0 & 0 & 0 & 0 \\
0 & 0 & 0 & 0 & (C(v_c^*))^{-1} & 0 & 0 & 0 & 0 \\
0 & 0 & 0 & 0 & 0 & (L(i_l^*))^{-1} & 0 & 0 & 0 \\
0 &  B_c &  0   &  0  & 0 & B_l & B_r R(i_r^*)  & B_j& 0 \\ 
0 &   0  &  Q_l & Q_m & Q_c & 0 & Q_r  & 0 & Q_u\end{pmatrix}
\end{eqnarray*}
and, via Schur's lemma, ${\cork F'(y^*, z^*)=1}$ follows again from Lemmas \ref{lema-siap} and
\ref{lema-loopscutsets}, which imply that
\begin{eqnarray*}
\ke \begin{pmatrix}
 B_c &  0   &   B_r R(i_r^*)  & B_j& 0 \\ 
   0  &  Q_l &  Q_r  & 0 & Q_u
\end{pmatrix}=\{0\}.
\end{eqnarray*}
As indicated above, this implies that the double zero eigenvalue is indeed index two.

\vs

\noindent 7. The fact that all non-vanishing eigenvalues 
($\lambda \neq 0$) have non-zero real part follows from the eigenvalue-eigenvector equations of
the pencil, which can be written as
\begin{eqnarray*} 
B_c u_c + \lambda B_l L(i_l^*)w_l + B_r R(i_r^*)w_r + B_j u_j  & = & 0 \\
 \lambda Q_c C(v_c^*) u_c + Q_l w_l + Q_r w_r +Q_m w_m + Q_u w_u  & = & 0 
\end{eqnarray*}
together with 
$\lambda \sigma_m =  w_m$, $w_c = \lambda C(v_c^*) u_c$, $u_l = \lambda L(i_l^*)w_l$. By taking
conjugate transposes and using the orthogonality of the cut and cycle spaces, we derive
\begin{equation}
\label{equPIEs}
 (\re \lambda) \left(u_c^{\star} C(v_c^*) u_c + w_l^{\star} L(i_l^*) w_l\right) + w_r^{\star}\frac{R(i_r^*)+(R(i_r^*))^{\tr}}{2}w_r =0.
\end{equation}

Now, if $\re \lambda = 0$, the positive definiteness of $R(i_r^*)$ implies $w_r=0$, 
and then
\begin{subequations}
\begin{eqnarray}\label{auxeq1}
 B_c u_c + \lambda B_l L w_l + B_j u_j & = & 0 \\
 \lambda Q_c C u_c + Q_l w_l + Q_m w_m + Q_u w_u & = & 0.\label{auxeq2}
\end{eqnarray}
\end{subequations}
But in this setting the hypothesis that there are no ILC-cutsets implies, in light of (\ref{auxeq1}), 
$u_c=0$ (and then $w_c=0$), 
$w_l=0$ (so that $u_l=0$)
and $u_j=0$. Additionally, the absence of VM-loops and (\ref{auxeq2}) would then imply
$w_m=0$
(and then $\sigma_m=0$) and $w_u=0$; this would yield a vanishing eigenvector, which is a
contradiction in terms.

\vs

\noindent 8. In order to check that 
condition 3 of Theorem \ref{th-daes} holds, we take $\oq$
from the requirement $F'(y^*, z^*)\oq \in \ke F'(y^*, z^*)-\{0\}$, which gives $\oq$ the form
$$\oq = ( \oq_1, \hspace{0.5mm} 0, \hspace{0.5mm} \oq_3, \hspace{0.5mm} \oq_4, \hspace{0.5mm} 0, 
\hspace{0.5mm} 0, \hspace{0.5mm} 0, \hspace{0.5mm} 0, \hspace{0.5mm} \oq_9), \text{  with } (\oq_3, 
\hspace{0.5mm} \oq_4, \hspace{0.5mm} \oq_9) \in \ke 
\begin{pmatrix} Q_l & \hspace{-1.5mm} Q_m & \hspace{-1.5mm}  Q_u \end{pmatrix}-\{0\}$$
and, in particular, ${\oq_4 \neq 0}$. The condition
$F''(y^*, z^*)\op \hspace{0.1mm}\oq \notin \im F'(y^*, z^*)$
then reads as
\begin{eqnarray*}
\hspace{-7mm}\begin{pmatrix} 0 \\ 0 \\  0 \\ B_m M'(q^*) \oq_4 \\ 0 \end{pmatrix} \notin 
\im \begin{pmatrix}
0  & 0  & 0  & 1  & 0  & 0  & 0  & 0  & 0 \\
0  & 0  & 0  & 0  & (C(v_c^*))^{-1}  & 0  & 0  & 0  & 0 \\
0  & 0  & 0  & 0  & 0  & (L(i_l^*))^{-1}  & 0  & 0  & 0 \\
0  &  B_c  &  0    &  0   & 0  & B_l  & B_r R(i_r^*)   & B_j & 0 \\ 
0  &   0   &  Q_l  & Q_m  & Q_c  & 0  & Q_r   & 0  & Q_u
\end{pmatrix}.
\end{eqnarray*}
Again, assuming that this condition is not met, we would derive the existence of a non-trivial
solution for
\begin{eqnarray*}
B_c u_2 -B_m M'(q^*)\oq_4 + B_r R(i_r^*) u_7 + B_i u_8 & = & 0,
\end{eqnarray*}
but the presence of a VML-loop with the memristor rules out,
because of the colored branch theorem 
(according to which, in 
a three-color graph with just one blue branch, this branch 
cannot form simultaneously a loop 
exclusively with green branches and a cutset only 
with red branches; cf.\ \cite{minty60, coloredbranch}),
the existence of any CMRI-cutset including the memristor. This would yield 
$M'(q^*) \oq_4 = 0$, against the fact that 
$M'(q^*) \neq 0  \neq \oq_4.$ 

\vs

\noindent 9. Finally, the proof that all eigenvalues (but the null one) have
negative real part at points of the equilibrium half-line where $M(q_m)>0$ 
(and, certainly, on a neighborhood of the bifurcating point) is essentially
similar to the one in item 7. Note only that the fact that $M(q_m) \neq 0$ make the 
eigenvalue-eigenvector equations read
\begin{eqnarray*} 
B_c u_c + \lambda B_l L(i_l)w_l + B_r R(i_r)w_r+ B_m M(q_m) w_m  + B_j u_j  & = & 0 \\
 \lambda Q_c C(v_c) u_c + Q_l w_l + Q_r w_r +Q_m w_m + Q_u w_u  & = & 0,
\end{eqnarray*}
and (\ref{equPIEs}) is now
\begin{equation*}
 (\re \lambda) \left(u_c^{\star} C(v_c) u_c + w_l^{\star} L(i_l) w_l\right) 
 + w_r^{\star}\frac{R(i_r)+(R(i_r))^{\tr}}{2}w_r+ w_m^{\star} M(q_m) w_m =0.
\end{equation*}
But since not only $C(v_c)$, $L(i_l)$, $R(i_r)$ but also $M(q_m)$ are positive 
definite (always at equilibrium points close enough to the bifurcating one), 
from the assumption $\re \lambda \geq 0$ we would derive $w_m = w_r = 0$ and
the reasoning proceeds as in item 7 above to show that all non-vanishing eigenvalues
must verify $\re \lambda < 0$. On the other hand, the fact that only one real eigenvalue
changes sign (and hence becomes positive) in the transition to the region where $M(q_m)<0$
follows from the TBWP phenomenon itself. This completes the proof.

\hfill $\Box$

\

\noindent Theorem \ref{th-circuits} shows that, in essence, the presence of a VML-loop (in
particular of an ML-loop) with the memristor and at least one inductor is the essential
configuration yielding a transcritical bifurcation without parameters, which occurs
if the memristance $M(q_m)$ vanishes (and eventually becomes negative) at a given $q_m^*$.
A simple parallel connection of a memristor and a linear inductor yields this
phenomenon, as shown in \cite{aims2014}.

Note that Theorem \ref{th-circuits} assumes that, except for the memristor, 
the remaining circuit devices are strictly locally passive. If this assumption is relaxed
things become more complicated; we present some results in this context in
the forthcoming section.

\section{Non-passive problems}
\label{sec-nonpassive}

\begin{figure}[htb]
\begin{center}
\epsfig{figure=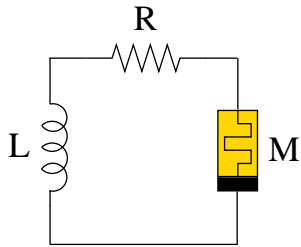, width=0.25\textwidth}
\vspace{-3mm}
\end{center}
\caption{MRL-circuit.}
\label{fig-MRL}
\end{figure}

Consider the circuit displayed in Figure \ref{fig-MRL}. Note that the absence
of a (V)ML-loop rules out an application of Theorem \ref{th-circuits}
in order to characterize an eventual transcritical bifurcation without
parameters in this circuit. However, it is easy to check that
the series connection of the linear resistor $R$ and the memristor $M$
can itself be modelled as a memristor with memristance $R+M(q)$. Moreover,
provided that $R+M(q^*)=0$ at a given $q^*$ (with $M'(q^*) \neq 0$), the 
circuit is expected and 
can be easily shown to undergo a TBWP.

\vv

Obviously, for the condition $R+M(q^*)=0$ to hold, either both $R$ and $M$ must
vanish or one of them need become negative. In particular, when $R$ becomes 
zero  or negative, the problem does not fit the strict passivity assumption
for resistors stated in Theorem \ref{th-circuits}. This means that a wider
framework is needed to address this phenomenon in general. 
Although
in its broad generality this is a difficult problem, some remarks in this
direction can be provided, as detailed in what follows. 

Specifically, we provide below 
conditions guaranteeing, in non-passive circuits with a single memristor,
that the null eigenvalue is indeed a multiple one, generalizing (as detailed
later, cf.\ the Remark after Proposition \ref{propo-nonpassive}) the framework considered in Theorem \ref{th-circuits}. We use the notion
of both a proper tree and an L-proper tree. Split
the branches of a given connected graph ${\cal G}$ in three pairwise disjoint sets $K_1$, $K_2$ and $K_3$,
in a way such that $K_1$ includes no loops and $K_3$ no cutsets. Then, 
as a consequence of the matroid structure of the set of acyclic subgraphs of ${\cal G}$
\cite{oxley}, one can guarantee that there exists at least one spanning tree including
all branches from $K_1$ and none from $K_3$ (an explicit proof can be found e.g.\ in
\cite{brown63}). Such a spanning tree is called (in general) a
{\em proper tree}. In circuit theory, this term is usually restricted to connected circuits without
VC-loops and IL-cutsets, to denote a spanning tree including all voltage sources and capacitors,
and neither current sources nor inductors. 
This notion can be traced back at least to \cite{bashkow57}. We will also make use of the
(in a certain sense dual) concept of an {\em L-proper tree}, which is a spanning tree 
including all voltage sources and inductors,
and neither current sources nor capacitors; such a tree exists if and only if the circuit
has neither  VL-loops nor IC-cutsets.
From \cite{cyclic} we borrow the concept of an MR-product; 
given a spanning tree in a connected
circuit with voltage and current sources, capacitors, inductors, resistors and memristors, 
the {\em MR-product} of this tree is simply the product of all resistances and memristances
in the co-tree branches (namely, the branches that do not belong to the spanning tree), evaluated
at equilibrium and
setting this product to 1 if all resistors and memristors are
actually located in the tree.

\begin{propo}  \label{propo-nonpassive}
Consider, as in Theorem \ref{th-circuits}, a circuit with a single memristor 
displaying an equilibrium
point at a given $(y^*,z^*)$, with 
$y^*=(q_m^*, v_c^*, i_l^*)$ and $z^*=(i_m^*, i_c^*, v_l^*, i_r^*, v_j^*, 
i_u^*)$. Assume that $C(v_c^*)$ and $L(i_l^*)$ are non-singular, besides
the following.
\begin{enumerate}
\item The circuit displays no VC-loops, IL-cutsets, VL-loops or IC-cutsets.
\item The sum of MR-products in proper trees does not vanish.
\item The sum of MR-products in L-proper trees does vanish.
\end{enumerate}
Then the algebraic multiplicity of the zero eigenvalue at $(y^*,z^*)$ is greater than one.
\end{propo}

The proof of this result follows from the results detailed in \cite{cyclic}.
First, the absence of VC-loops and IL-cutsets,
together with the non-vanishing condition on the MR-product sum in proper trees, guarantees the
matrix $g_z(y^*, z^*)$ to be non-singular, cf.\ Proposition 3 in \cite{cyclic}. 
Note that this matrix has now the form
\begin{eqnarray}
\label{matr1}
g_z(y^*, z^*)=\begin{pmatrix}B_m M(q_m^*) & 0 & B_l & B_r R(i_r^*)  & B_j& 0 \\ Q_m & Q_c & 0 & Q_r  & 0 & Q_u\end{pmatrix}.
\end{eqnarray}
Additionally, the null eigenvalue having a multiplicity greater than one 
is in a way the dual property to the one above, and
relies on the structure of 
the matrix $K(0)$ from (\ref{bigmatrix}), which in the presence of a possibly non-vanishing memristance
reads as
\begin{eqnarray*}
\begin{pmatrix}
0 & 0 & 0 & -(C(v_c^*))^{-1} & 0 & 0 & 0 & 0 \\
 0 & 0 & 0 & 0 & -(L(i_l^*))^{-1} & 0 & 0 & 0 \\
 -B_c &  0   &  -B_m M(q_m^*) & 0 & -B_l & -B_r R(i_r^*)  & -B_j& 0 \\ 
   0  &  -Q_l & -Q_m & -Q_c & 0 & -Q_r  & 0 & -Q_u
\end{pmatrix}.
\end{eqnarray*}
Again, this matrix can be checked
to be singular if and only if so it is the matrix
\begin{eqnarray}
\begin{pmatrix}
 B_c &  0   &  B_mM(q_m^*)  &  B_r R(i_r^*)  & B_j& 0 \\
   0  &  Q_l & Q_m &  Q_r  & 0 & Q_u
\end{pmatrix},  \label{matr2}
\end{eqnarray}
but from Proposition 3 in \cite{cyclic} one can show that,
in the absence of VL-loops and IC-cutsets, the latter matrix is singular if and
only if the sum of MR-products in L-proper trees does vanish. 
Details are not difficult and 
are left to the reader. The above-referred duality property becomes clear in the light of 
the matrices
(\ref{matr1}) and (\ref{matr2}), which have exactly the same form after a change of reactive
devices (capacitors and inductors) and an obvious column reordering.

Proposition \ref{propo-nonpassive} opens a way for future research, which should address 
the remaining conditions from Theorem \ref{th-daes} in order
to characterize the TBWP in this
wider setting. Note that Proposition \ref{propo-nonpassive}
does not require any passivity assumption on the circuit characteristic matrices. 

\vspace{2mm}

\noindent  {\bf Remark.} In the presence of an VML-loop, 
as in the setting of Theorem \ref{th-circuits},
the memristor must by definition belong to the cotree of all L-proper trees (since these must accommodate
all voltage sources and inductors); $M$ is therefore a common factor in all MR-products and
therefore the condition $M=0$ in Theorem \ref{th-circuits} arises naturally.

\vspace{2mm}

\begin{figure}[h]
\vspace{6mm}
\centering
\ \mbox{} \hspace{22mm} \
\parbox{2.8in}{
\epsfig{figure=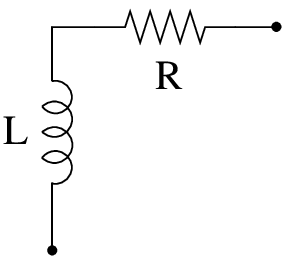, width=0.17\textwidth}
}
\ \hspace{-27mm} \parbox{2.8in}{\vspace{1.5mm}
\epsfig{figure=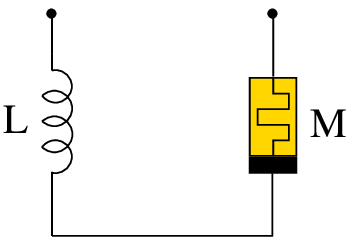, width=0.215\textwidth}
}\ \hspace{-15mm}
\vspace{2mm}
\caption{L-proper trees in the circuit of Figure \ref{fig-MRL}}
\label{fig-properMR}
\vspace{1mm}
\end{figure}

Proposition \ref{propo-nonpassive} above explains, in topological terms, why the
bifurcating condition in the circuit of Figure \ref{fig-MRL} is $M(q_m^*)+R = 0$. Just
note that the circuit has two L-proper trees, displayed in Figure \ref{fig-properMR},
the cotrees of which amount, respectively, to the memristor $M$ and the resistor $R$;
the sum of products arising in item 3 of Proposition \ref{propo-nonpassive} is just $M(q_m^*) + R$
and for this reason the vanishing of this sum is responsible for a multiple zero eigenvalue 
supporting the TBWP.
However, the scope of this results goes much further than this elementary (pedagogic) example, since
the computation of spanning trees is an easily automatable task and therefore 
the result applies to much more complex circuits. A case of intermediate complexity is discussed
below with illustrative purposes.

\subsection*{A memristive artificial neural network}

Memristors provide an excellent framework for the implementation of artificial
neural networks. The key reason is that they define a nanometric scale, electrically
adaptable device perfectly suited for implementing neural synapses, emulating (to a certain
extent) the STDP mechanism from biological neural systems \cite{serrano2013}.
A lot of recent
literature explores this idea; see e.g.\
\cite{adhikari2012, jo2010, kim2012bis, pershin10c, soudry, wen} and references therein.

\begin{figure}[h]
\ \hspace{4mm}
\parbox{3in}{\hspace{0mm}
\epsfig{figure=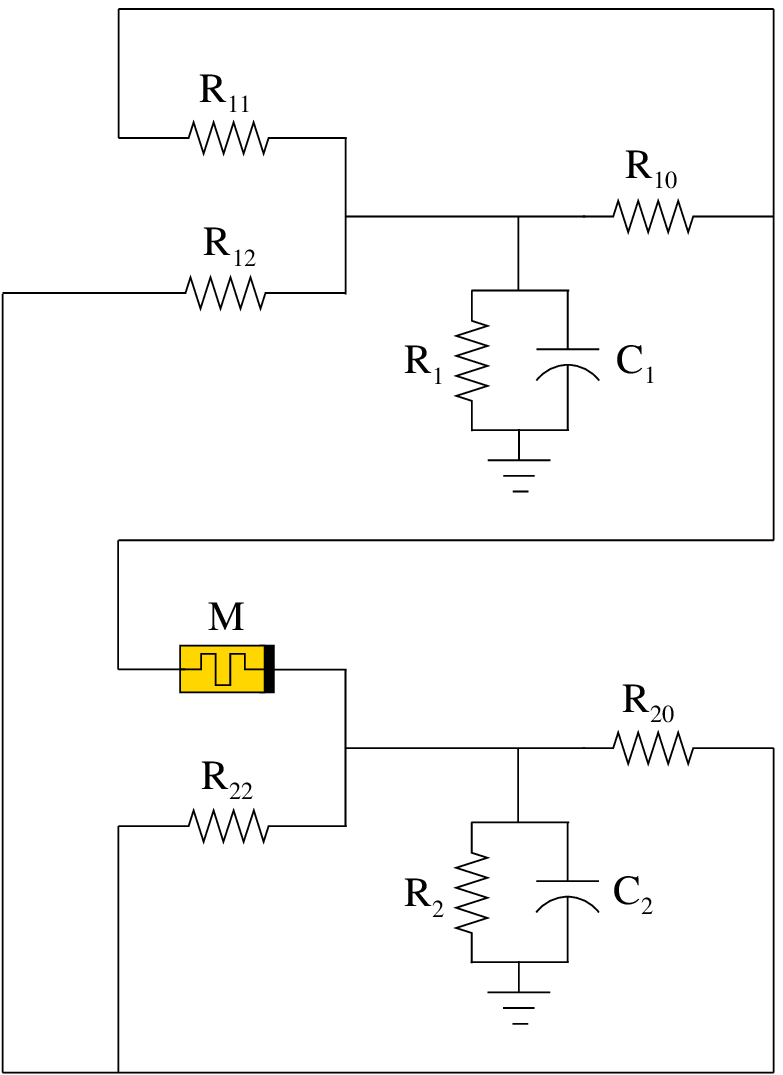, width=0.40\textwidth}
}
\hspace{2mm}
\parbox{3in}{ \vspace{41.5mm}
\epsfig{figure=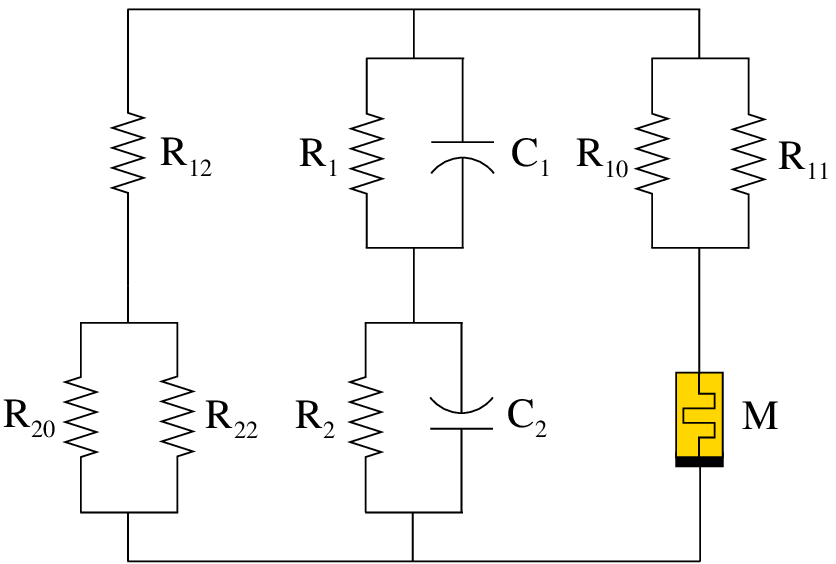, width=0.44\textwidth}
}
\vspace{2mm}
\caption{(a) Memristive network. \hspace{15mm} (b) Alternative circuit description.}
\label{fig-neural}
\end{figure}

In this context, we analyze below a simplification of an 
additive model proposed in \cite{wen}. We ignore delays and assume that 
each neuron is defined by a passive RC-connection, and also that
the input-output function in each neuron is just implemented by a linear
passive resistor. Furthermore, we focus on a problem with 
just two neurons and assume
that the conductivities of three out of four synaptic connections are fixed,
in order to concentrate the attention on a bifurcating memristor.
This simplified
model is depicted in Figure \ref{fig-neural} (a); in (b) 
we arrange the circuit in a more convenient manner for later computations.

\vvvvv

Note that, for any charge $q_m$ in the memristor and null values
in the remaining circuit variables one gets an equilibrium point.
This simply expresses the presence of a line
of equilibria along the $q_m$-coordinate axis,
and the linearized dynamics exhibits a zero eigenvalue
along this axis. 
Omitting details for the sake of brevity, a state-space model
for this circuit indicates that 
the condition for the zero eigenvalue to be multiple is
\begin{eqnarray} \label{bifcond}
R_A R_B + (R_1 + R_2)\left[(R_{20} +R_{22})R_A
+ (R_{10} + R_{11})R_B\right] = 0,
\end{eqnarray}
with
\begin{eqnarray*}
R_A = M(R_{10}+ R_{11}) + R_{10}R_{11},
\ R_B = R_{12}(R_{20} +  R_{22}) + R_{20}R_{22}.
\end{eqnarray*}
This yields the bifurcation value
\begin{eqnarray*}
M = \frac{- R_C -R_D}{(R_{10}+ R_{11})[R_B +(R_1 + R_2)(R_{20} +R_{22})]}
\end{eqnarray*}
where
\begin{eqnarray*}
R_C = R_{10}R_{11}[R_B +(R_1 + R_2)(R_{20} +R_{22})], 
\ R_D = R_B(R_1 + R_2)(R_{10} + R_{11}).
\end{eqnarray*}
Actually, for negative values $M$ smaller than the one above, and provided that all resistances are positive,
a transition of an eigenvalue to the positive real semiaxis 
signals a stability loss due to to a TBWP.



\begin{figure}[ht]
\parbox{0.7in}{
\epsfig{figure=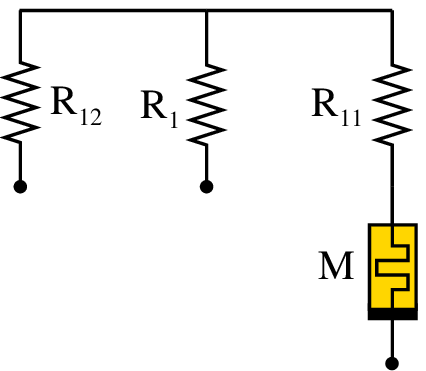, width=0.18\textwidth}
}
\hspace{12mm}
\parbox{0.7in}{
\epsfig{figure=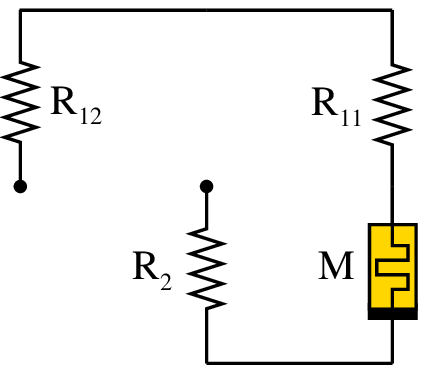, width=0.18\textwidth}
}
\hspace{12mm}
\parbox{0.7in}{
\epsfig{figure=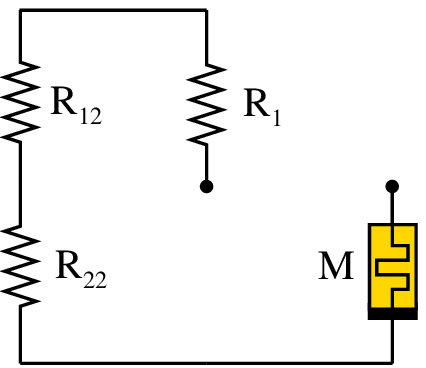, width=0.18\textwidth}
}
\hspace{12mm}
\parbox{0.7in}{
\epsfig{figure=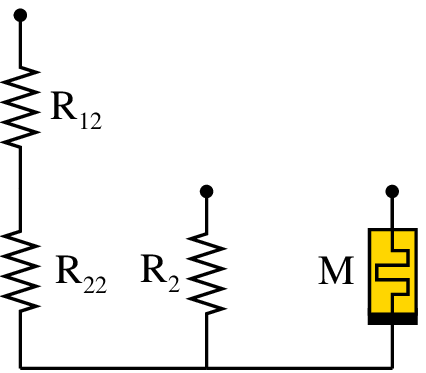, width=0.18\textwidth}
}
\hspace{12mm}
\parbox{0.7in}{
\epsfig{figure=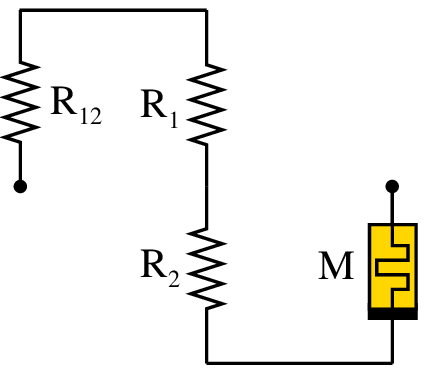, width=0.18\textwidth}
}\vspace{5mm}\\
\parbox{0.7in}{
\epsfig{figure=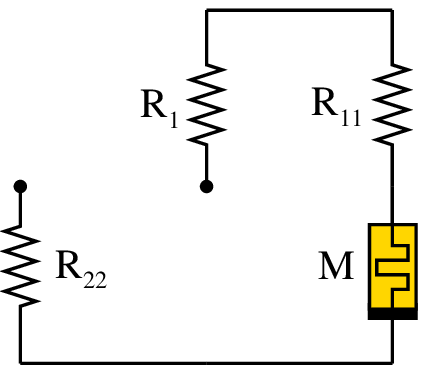, width=0.18\textwidth}
}
\hspace{12mm}
\parbox{0.7in}{
\epsfig{figure=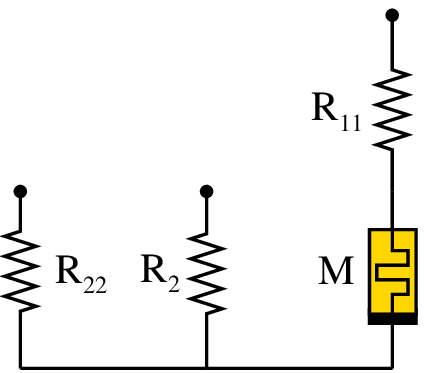, width=0.18\textwidth}
} 
\hspace{12mm}
\parbox{0.7in}{
\epsfig{figure=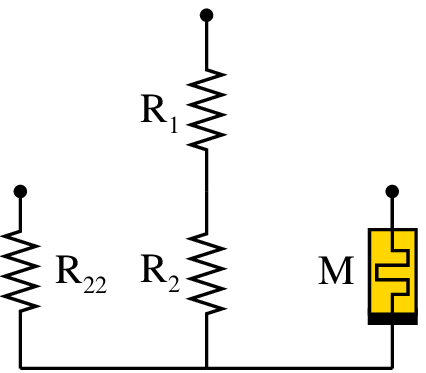, width=0.18\textwidth}
}
\hspace{12mm}
\parbox{0.7in}{
\epsfig{figure=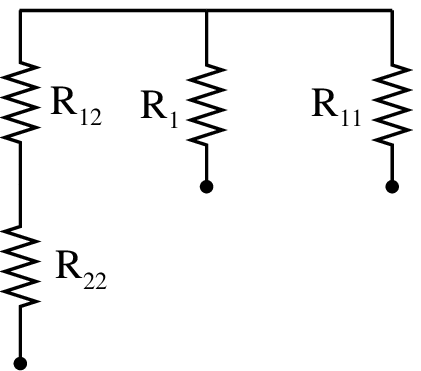, width=0.18\textwidth}
}
\hspace{12mm}
\parbox{0.7in}{
\epsfig{figure=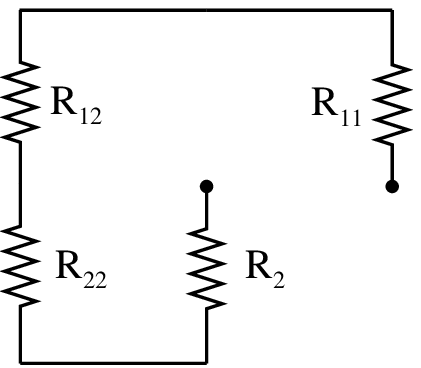, width=0.18\textwidth}
}\vspace{5mm}\\
\parbox{0.7in}{
\epsfig{figure=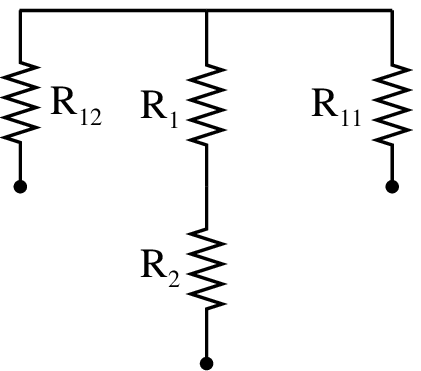, width=0.18\textwidth}
}
\hspace{12mm}
\parbox{0.7in}{
\epsfig{figure=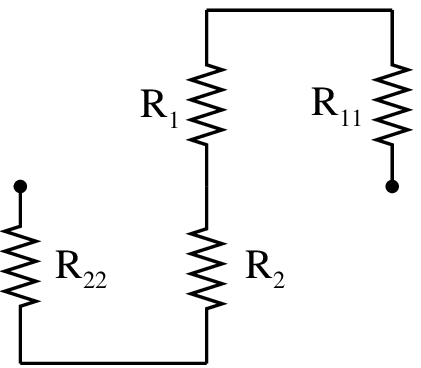, width=0.18\textwidth}
}
\hspace{12mm}
\parbox{0.7in}{
\epsfig{figure=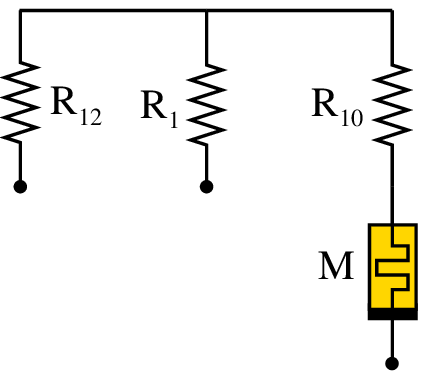, width=0.18\textwidth}
}
\hspace{12mm}
\parbox{0.7in}{
\epsfig{figure=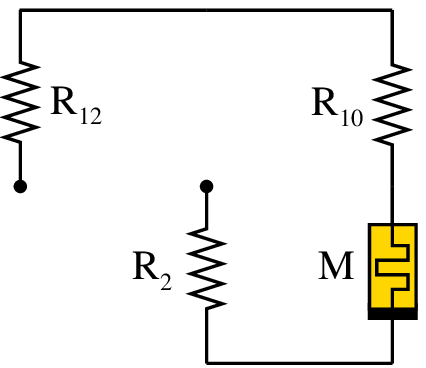, width=0.18\textwidth}
}
\hspace{12mm}
\parbox{0.7in}{
\epsfig{figure=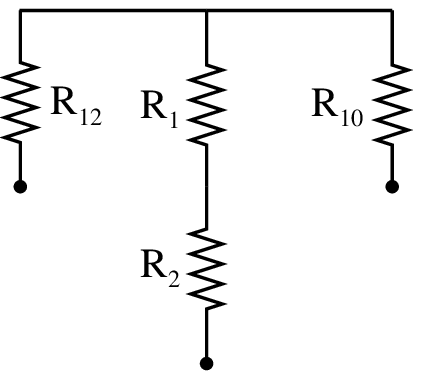, width=0.18\textwidth}
}\vspace{5mm}\\
\parbox{0.7in}{
\epsfig{figure=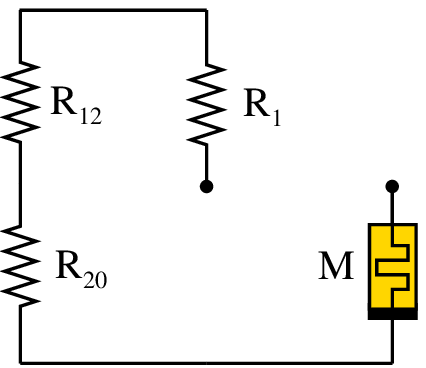, width=0.18\textwidth}
}
\hspace{12mm}
\parbox{0.7in}{
\epsfig{figure=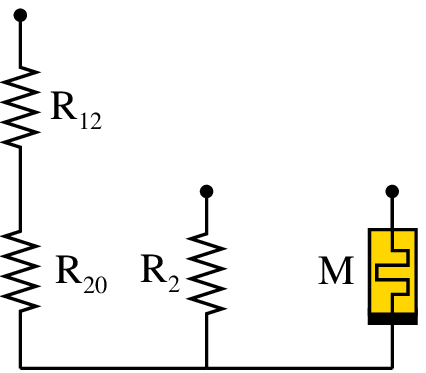, width=0.18\textwidth}
} 
\hspace{12mm}
\parbox{0.7in}{
\epsfig{figure=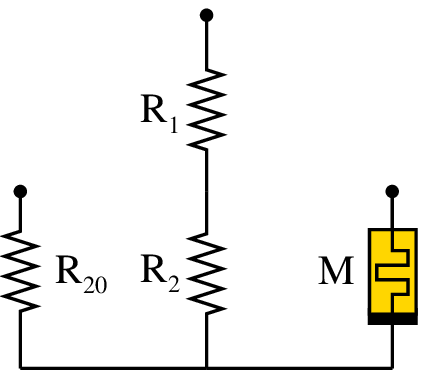, width=0.18\textwidth}
}
\hspace{12mm}
\parbox{0.7in}{
\epsfig{figure=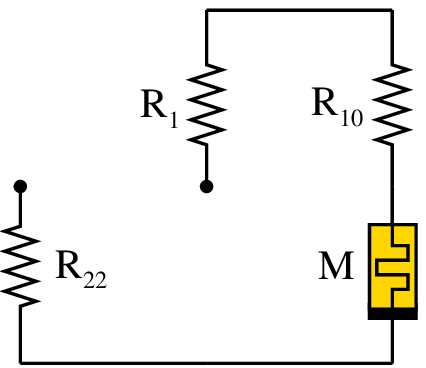, width=0.18\textwidth}
}
\hspace{12mm}
\parbox{0.7in}{
\epsfig{figure=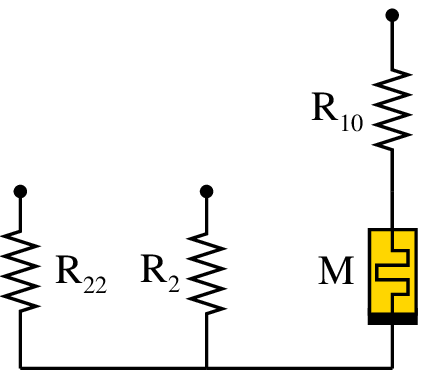, width=0.18\textwidth}
}\vspace{5mm}\\
\parbox{0.7in}{
\epsfig{figure=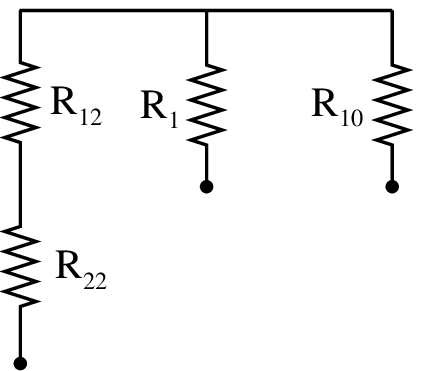, width=0.18\textwidth}
}
\hspace{12mm}
\parbox{0.7in}{
\epsfig{figure=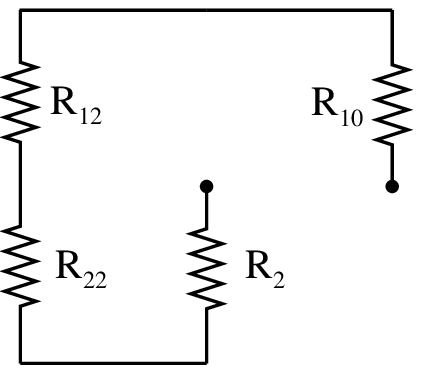, width=0.18\textwidth}
}
\hspace{12mm}
\parbox{0.7in}{
\epsfig{figure=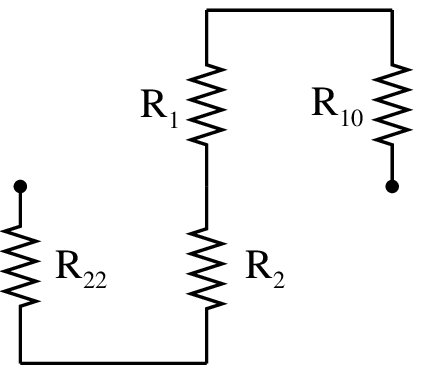, width=0.18\textwidth}
}
\hspace{12mm}
\parbox{0.7in}{
\epsfig{figure=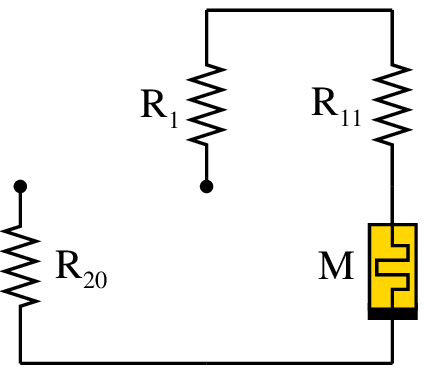, width=0.18\textwidth}
}
\hspace{12mm}
\parbox{0.7in}{
\epsfig{figure=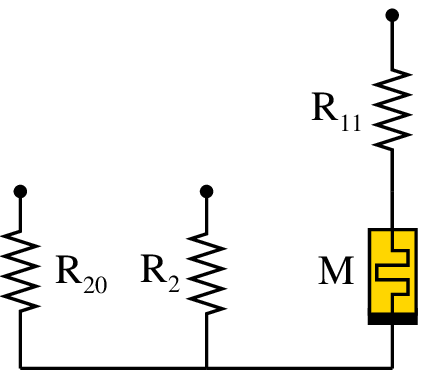, width=0.18\textwidth}
}\vspace{5mm}\\
\parbox{0.7in}{
\epsfig{figure=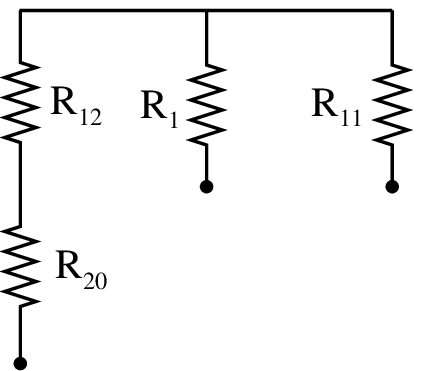, width=0.18\textwidth}
}
\hspace{12mm}
\parbox{0.7in}{
\epsfig{figure=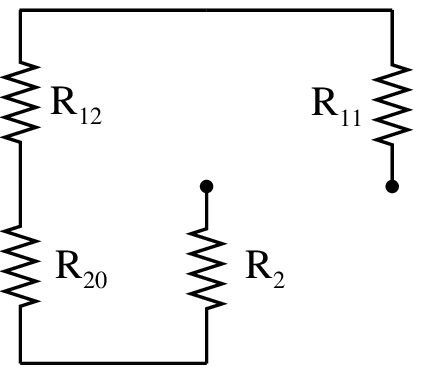, width=0.18\textwidth}
} 
\hspace{12mm}
\parbox{0.7in}{
\epsfig{figure=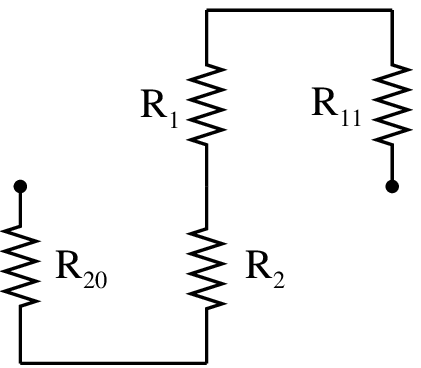, width=0.18\textwidth}
}
\hspace{12mm}
\parbox{0.7in}{
\epsfig{figure=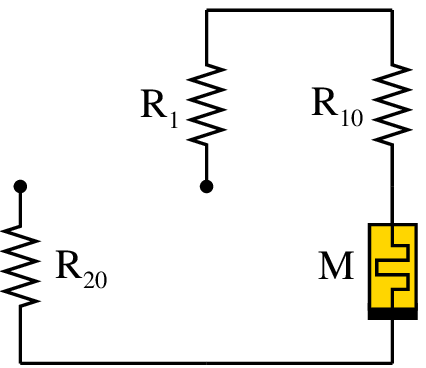, width=0.18\textwidth}
}
\hspace{12mm}
\parbox{0.7in}{
\epsfig{figure=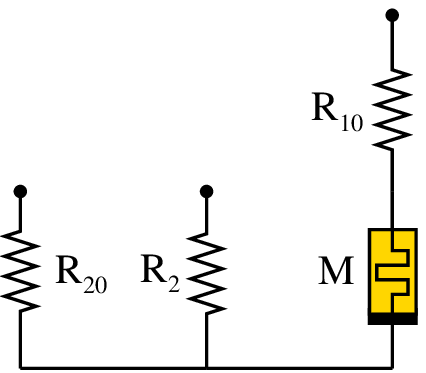, width=0.18\textwidth}
}\vspace{5mm}\\
\parbox{0.7in}{
\epsfig{figure=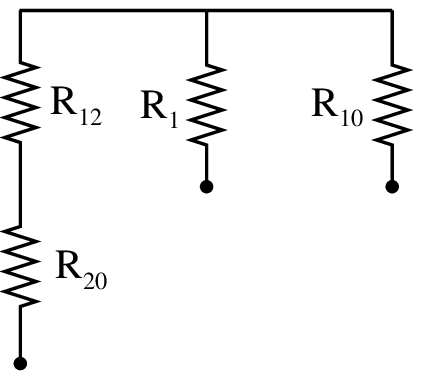, width=0.18\textwidth}
}
\hspace{12mm}
\parbox{0.7in}{
\epsfig{figure=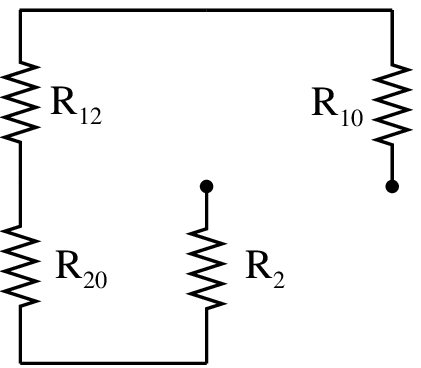, width=0.18\textwidth}
}
\hspace{12mm}
\parbox{0.7in}{
\epsfig{figure=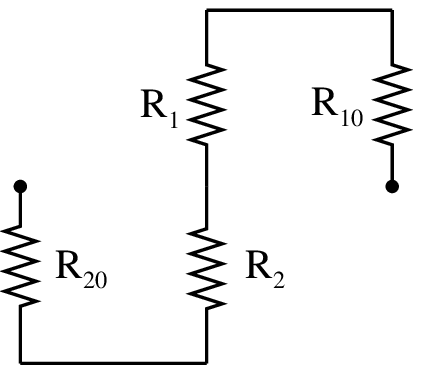, width=0.18\textwidth}
}\vspace{5mm}
\caption{L-proper trees for the circuit in Figure \ref{fig-neural}.} \label{fig-trees}
\end{figure}

Our goal is to 
explain in topological terms the condition (\ref{bifcond})
on the resistances and the memristance, 
which makes this null eigenvalue
a multiple one,
in order to illustrate the scope of Proposition \ref{propo-nonpassive}. 
Note, that in large scale circuits the derivation of a 
model and therefore the explicit computation
of the bifurcation conditions are usually unfeasible, 
and for this reason one has no option but to resort to circuit-theoretic 
results such as the one
in Proposition \ref{propo-nonpassive}.
To achieve this one needs to compute the set of L-proper trees
of the circuit. In our present example there are actually 
33 L-proper trees, depicted in Figure \ref{fig-trees}. 


\vspace{2mm}

One can check that this set of trees actually explains the
bifurcating condition (\ref{bifcond}), in light of item 3
of Proposition \ref{propo-nonpassive}. Specifically,
the terms responsible for the 
product $R_AR_B$ arise from the cotree branches of 
the trees 5, 8, 11, 12, 15, 18, 23, 28 and 33 in Figure \ref{fig-trees}; analogously, the terms
$(R_1 + R_2)(R_{20} +R_{22})R_A$ come from the cotree branches
of trees 3, 4, 9, 10, 16, 17, 21, 22, 26, 27, 31 and 32, 
whereas the products $(R_1 + R_2)(R_{10} + R_{11})R_B$ arise
from the cotrees of the trees 
1, 2, 6, 7, 13, 14, 19, 20, 24, 25, 29 and 30. 

As indicated
above, this graph-theoretic characterization of the bifurcation
condition is easily scalable to
large circuits in which the analytical (model-based) computation
of such condition is not feasible.

\section{Concluding remarks}
\label{sec-con}

We have presented in this paper a detailed circuit-theoretic
characterization of the transcritical bifurcation without parameters
in circuits with one memristor, systematically yielding
lines of equilibrium points. To do so, we have developed 
mathematical statements of the TBWP theorem for explicit ODEs 
in arbitrary dimension and also for semiexplicit DAEs, which are believed to
be of independent interest. This allows
for a graph-theoretic analysis of the bifurcation in the circuit context.
Future research should provide a complete characterization of 
this phenomenon in non-passive settings, along the lines
discussed in Section \ref{sec-nonpassive}. Other related
bifurcations, such as the Hopf bifurcation without parameters,
might be analyzed in similar terms.

\section*{Appendix: digraph matrices}

In the formulation of the circuit model (\ref{branch}) we
make use of the so-called {\em loop} and {\em cutset} matrices
defined below. 
Given a digraph with $m$ edges, $n$ nodes and $k$ connected
components, choose an orientation in every loop and define
componentwise the
{\em loop matrix} $\tilde{B}$ as $(b_{ij})$, with
\begin{eqnarray*}
b_{ij} = \left\{
\begin{array}{rl}
1 & \text{ if edge } j \text{ is in loop } i \text{ with the same
  orientation }\\
-1 & \text{ if edge } j \text{ is in loop } i \text{ with the opposite
  orientation } \\
0 & \text{ if edge } j \text{ is not in loop } i.
\end{array}
\right.
\end{eqnarray*}
This matrix has rank
$m-n+k$, and a {\em reduced loop matrix} $B$
is any $((m-n+k)\times m)$-submatrix of $\tilde{B}$ with full row rank.

The dual concept is that of a reduced cutset
matrix. Recall that a set $K$ of edges in a 
digraph is a {\em cutset} if the 
removal of $K$ increases the number of connected components,
and $K$ is minimal with respect to this property, that is,
retaining one or more edges from $K$ keeps the number of
components invariant. 
All the edges of a cutset may be shown to connect the same
pair of connected components which result from the cutset deletion, 
and this allows one to define the orientation of a cutset, say
from one of these components towards the other.
This makes it possible to define the cutset matrix $\tilde{Q} = (q_{ij})$ as
\begin{eqnarray*}
q_{ij} = \left\{
\begin{array}{rl}
1 & \text{ if edge } j \text{ is in cutset } i \text{ with the same
  orientation }\\
-1 & \text{ if edge } j \text{ is in cutset } i \text{ with the opposite
  orientation } \\
0 & \text{ if edge } j \text{ is not in cutset } i.
\end{array}
\right.
\end{eqnarray*}
Now the rank of $\tilde{Q}$ is $n-k$, and 
a {\em reduced cutset matrix} $Q \in \R^{(n-k) \times m}$ is obtained
by choosing  any set of 
$n-k$ linearly independent rows of $\tilde{Q}$.

\end{document}